\documentclass[11pt]{article}
\usepackage{amsfonts}
\usepackage{mathrsfs,color,amscd,amssymb,enumerate,amsthm,amsmath,bm,graphicx,psfrag,subfigure}
\usepackage{latexsym}
\usepackage{float,fancybox,shapepar,setspace,hyperref}
\usepackage{pgf,tikz}

\makeatletter
\def\leftharpoonfill@{\arrowfill@\leftharpoonup\relbar\relbar}
\def\rightharpoonfill@{\arrowfill@\relbar\relbar\rightharpoonup}
\newcommand\rbjt{\mathpalette{\overarrow@\rightharpoonfill@}}
\newcommand\lbjt{\mathpalette{\overarrow@\leftharpoonfill@}}
\makeatother

\makeatletter

\renewcommand{\@seccntformat}[1]{{\csname the#1\endcsname}{\normalsize .}\hspace{.5em}}
\makeatother

\def \[{\begin{equation}}
\def \]{\end{equation}}

\def\qed{ \hfill $\square$}
\def\bqed{ \hfill $\blacksquare$}

\newtheorem{thm}{Theorem}
\newtheorem{prop}{Proposition}

\newtheorem{fac}{Fact}
\newtheorem{lem}{Lemma}

\newtheorem{conj}{Conjecture}

\usetikzlibrary{arrows}
\begin{document}

\title{Improved  bound on the number of edges of diameter-$k$-critical graphs}
\author{Xiaolin Wang\thanks{School of Mathematics and Statistics, Fuzhou University, Fuzhou
350108, P.R. China (xiaolinw@fzu.edu.cn)}, Yanbo Zhang\thanks{ Corresponding author. School of Mathematical Sciences, Hebei Normal University, Shijiazhuang
050024, P.R. China (Email: ybzhang@hebtu.edu.cn)}~, Xiutao Zhu\thanks{School of Mathematics, Nanjing University of Aeronautics and Astronautics, Nanjing 211106, P.R.
China (Email: zhuxiutao@nuaa.edu.cn)}}

\date{}
\maketitle

\begin{abstract}
A graph  is  diameter-$k$-critical  if its diameter equals   $k$ and the deletion of any edge increases its diameter. The
Murty-Simon Conjecture states that
for any diameter-2-critical graph $G$ of order $n$, $e(G) \leq \lfloor \frac{n^2}{4}\rfloor$,  with equality   if and only if $G \cong K_{\lfloor \frac{n}{2}\rfloor,\lceil \frac{n}{2}\rceil}$.
F\"uredi
(JGT,1992) proved that this conjecture is true for sufficiently large $n$. Over two decades later, Loh and Ma (JCT-B, 2016) proved that $e(G) \leq \frac{n^2}{6}+o(n^2)$ for  diameter-3-critical graphs $G$, and
   $e(G) \leq \frac{3n^2}{k}$ for diameter-$k$-critical graphs $G$ with $k \geq 4$. In this paper, we improve the  bound for diameter-$k$-critical graphs  to   $ \frac{n^2}{2k}+o(n^2)$.

\vskip 2mm

\end{abstract}
{\setcounter{section}{0}

\section{Introduction}\setcounter{equation}{0}
\vskip 2mm
Let $G=(V(G),E(G))$ be a finite simple connected graph, and  $e(G)$ denote the number of edges in $G$, i.e., $e(G)=|E(G)|$.
A \emph{(u,v)-path} is a path with endpoints $u$ and $v$, and its \emph{length} is the number of its edges.
For any $u,v\in V(G)$, the distance $d_G(u,v)$ is the length of a shortest $(u,v)$-path  in $G$. If  no $(u,v)$-path exists, we say the distance $d_G(u,v)=\infty$. The  \textit{diameter} of $G$ is defined as $diam(G)=\max\{d(u,v)\,|\,u,v \in V(G)\}$. We say $G$ is \emph{diameter-k-critical}  if $diam(G)=k$ and $diam(G-e)>k$ for any $e\in E(G)$.

The research of the number of edges of diameter-$k$-critical graphs is one of the oldest subjects of the extremal graph theory. In the 1960s and 1970s, Ore \cite{O}, Plesn\'ik  \cite{P}, and Murty and Simon (see in \cite{CH}) independently presented the following long-standing conjecture, commonly known as the Murty-Simon Conjecture.
\begin{conj}\label{MS}
For any diameter-2-critical graph $G$ with $n$ vertices, $e(G) \leq \lfloor \frac{n^2}{4}\rfloor$. Moreover,  equality holds if and only if $G \cong K_{\lfloor \frac{n}{2}\rfloor,\lceil \frac{n}{2}\rceil}$.
\end{conj}

In 1979, Cacetta and H\"aggkvist \cite{CH} proved  that for any diameter-2-critical graph $G$, $e(G) \leq 0.27n^2$. In 1987, Fan \cite{F}  proved that if $n\leq 24$ or $n=26$, then $e(G) \leq 0.25n^2$; and if $n \geq 25$, then $e(G) \leq 0.2532n^2$.  In 1992, F\"uredi \cite{Fu} made a significant breakthrough, which involved a clever application of the Ruzsa-Sezemer\'edi $(6,3)$ theorem to prove that the Murty-Simon Conjecture is true   for large $n$ (nonasymptotic). Since the current quantitative bounds on the $(6,3)$ theorem are of tower type, the constraint of $n$ is quite intense. Hence, it is interesting to find another way that avoids Regularity-type ingredients to solve Conjecture \ref{MS} completely. In the 2000s and 2010s, Haynes, Henning, van der Merwe, Yeo, et al. studied  the diameter-2-critical graphs via its complement graphs, and proved that the conjecture is true under some constraints; see \cite{H1}-\cite{H8}.

For $k \geq 3$,
 Cacetta and H\"aggkvist \cite{CH} constructed a class of diameter-$k$-critical graphs  and conjectured that they have the maximum number of edges. However, in 1981, Krishnamoorthy and Nandakumar \cite{KN} presented a similar class of diameter-$k$-critical graphs  with larger number of edges,
and conjectured that these graphs maximize the number of edges: for any $k \geq 3$,
  take $a_1$ disjoint paths $x_1^ix_2^i \cdots x_{k-1}^i$ for $i=1,2,\ldots,a_1$, join each $x_1^i$ to $a_0$ new vertices, and  each $x_{k-1}^i$ to $a_2$ new vertices.
Loh and Ma \cite{LM} proved that these graphs are
diameter-$k$-critical for any positive integers $a_0,a_1,a_2$.
 Krishnamoorthy and Nandakumar \cite{KN} observed that if $a_0=1$, $a_1\approx \frac{n}{2(k-1)}$ and $a_2=n-a_0-a_1(k-1)$, then these graphs achieve the maximum number of edges, i.e., $\frac{n^2+2(k-2)n}{4(k-1)}+O(1)$. In fact, it is easy to verify that  $a_1\approx \frac{n}{2(k-1)}$ is enough to guarantee the same number of edges. For convenience, denote these graphs with $a_1\approx \frac{n}{2(k-1)}$ as $G_k$ for all $k \geq 3$.

Surprisingly, for $k=3$,  $G_3$ is not the unique diameter-3-critical graph which attains the bound $\frac{n^2+2n}{8}+O(1)$. Loh and Ma \cite{LM} found another kind of diameter-3-critical graph $G_{3,0}$ with even order $n$. This graph is constructed by
adding a perfect matching between $K_{\frac{n}{2}}$ and its complement $\overline{K_{\frac{n}{2}}}$.
 It is easy to check that $G_{3,0}$ is diameter-3-critical and $e(G_{3,0})=e(G_3)=\frac{n^2+2n}{8}$. Since the structures of $G_3$ and  $G_{3,0}$ are quite different,
 the study of diameter-3-critical graphs may be more difficult than other diameter-$k$-critical graphs. In this paper, we  identify some  other diameter-3-critical graphs with  $\frac{n^2+2n}{8}$ edges.

\vskip 2mm
Loh and Ma \cite{LM} also  provided the bounds of the number of edges of  diameter-$k$-critical graphs.
They proved the following theorem.

\begin{thm}\label{n/6}
(Loh and Ma \cite{LM}) For any diameter-3-critical graph $G$ on $n$ vertices, $e(G) \leq \frac{n^2}{6}+o(n^2)$; and for any diameter-$k$-critical graph $G$ on $n$ vertices,
$e(G) \leq \frac{3n^2}{k}$.
\end{thm}

In this paper, we improve their bound, and the proof also gives the same result for diameter-3-critical graphs, though the approach is different from Loh and Ma's \cite{LM}.

\begin{thm}\label{n^2/2k}
Let $G$ be a diameter-$k$-critical graph on $n$ vertices with $k \geq 3$. Then we have $e(G) \leq \frac{n^2}{2k}+k^2o(n^2)$.
\end{thm}



The rest of this paper is organized as follows.
Section 2 presents another class of diameter-3-critical graphs with $\frac{n^2+2n}{8}$ edges. In Section 3, we prove two lemmas for graphs with diameter $k$. The proof of Theorem \ref{n^2/2k} is in Section 4.

\section{Diameter-3-critical graphs with $\frac{n^2+2n}{8}$   edges}\setcounter{equation}{0}\label{graphs}

In this section, we construct a class of graphs from $G_{3,0}$ which are diameter-3-critical with $\frac{n^2+2n}{8}$ edges.
 Let $n \geq 6$  be even, $\theta$ be a vertex bijection between $V(K_{\frac{n}{2}})$ and  $V(\overline{K_{\frac{n}{2}}})$ defined by the perfect matching in $G_{3,0}$.  Choose any matching $M$ in $K_{\frac{n}{2}}$, delete $M$ and add  each edge $\theta(u)\theta(v)$ to $\overline{K_{\frac{n}{2}}}$ if and only if $uv \in M$. We denote such a graph as $G_{3,M}$. Obviously, $e(G_{3,M})=e(G_{3,0})$ for any matching $M$ in $K_{\frac{n}{2}}$. We now prove these  graphs are  diameter-3-critical.
\begin{prop}
For any matching $M$ in $K_{\frac{n}{2}}$, $G_{3,M}$ is a diameter-3-critical graph.

\end{prop}

\noindent{\bf{Proof.}}
We first prove that $diam(G_{3,M})=3$. For any two vertices $u,v \in V(K_{\frac{n}{2}})$, if $uv \notin M$, then $d(u,v)=1$. If $uv \in M$, then there exists another vertex $w \in V(K_{\frac{n}{2}})$ such that $uwv$ is a path in
$G_{3,M}$, implying that $d(u,v)=2$.
For $u \in V(K_{\frac{n}{2}})$ and $\overline{v} \in V(\overline{K_{\frac{n}{2}}})$, if $\theta(u)=\overline{v}$, then $d(u,\overline{v})=1$. Otherwise, either $u\theta^{-1}(\overline{v})\overline{v}$ or
$u\theta(u)\overline{v}$
is a path in $G_{3,M}$, implying that $d(u,\overline{v})=2$.
For any two vertices $\overline{u},\overline{v} \in V(\overline{K_{\frac{n}{2}}})$,
if $\theta^{-1}(\overline{u})\theta^{-1}(\overline{v}) \in M$, then $d(\overline{u},\overline{v})=1$.
If $\theta^{-1}(\overline{u})\theta^{-1}(\overline{v}) \notin M$, then $\overline{u}\theta^{-1}(\overline{u})\theta^{-1}(\overline{v})\overline{v}$ is a path in $G_{3,M}$, implying that $d(\overline{u},\overline{v})\leq 3$. Since  $\overline{u}$ and $\overline{v}$ are nonadjacent and
have no common neighbor,  $d(\overline{u},\overline{v})= 3$. Note that there must exist
$\overline{u},\overline{v} \in V(\overline{K_{\frac{n}{2}}})$ such that
$\theta^{-1}(\overline{u})\theta^{-1}(\overline{v}) \notin M$, then
 $diam(G_{3,M})=3$.

Now we  prove its critical property. Note that $E(G_{3,M}[V(\overline{K_{\frac{n}{2}}})])$ is a matching.
For any edge $uv \in E(G_{3,M}[V(K_{\frac{n}{2}})])$,   there exists exactly one path
$\theta(u)uv\theta(v) $ of length at most 3 connecting $\theta(u)$ and $\theta(v)$, that is, $d_{G_{3,M}-uv}(\theta(u),\theta(v)) >3$.
For any edge $\overline{u}\overline{v} \in E(G_{3,M}(V(\overline{K_{\frac{n}{2}}})))$,  since $\theta^{-1}(\overline{u})\theta^{-1}(\overline{v}) \notin E(G_{3,M})$, there exists exactly one path
$\overline{u}\overline{v} $ of length at most 3 connecting  $\overline{u}$ and $\overline{v}$, that is, $d_{G_{3,M}-\overline{u}\overline{v}}(\overline{u},\overline{v}) >3$.
For any edge $u\theta(u)\in E(G_{3,M})$, $u \in V(K_{\frac{n}{2}})$ and $\theta(u) \in V(\overline{K_{\frac{n}{2}}})$, we can check that whether $u$ is  covered by $M$ or not, there exists exactly one path
$u\theta(u) $ of length at most 3 connecting  $u$ and $\theta(u)$, that is, $d_{G_{3,M}-u\theta(u)}(u,\theta(u)) >3$.
\bqed

\section{Two lemmas for  graphs with diameter $k$}\setcounter{equation}{0}\label{somelemmas}

In this section, we  generalize two  lemmas of F\"{u}redi \cite{Fu}, and of  Loh and Ma \cite{LM} to any graph $G$ with diameter $k \geq 2$ (no need to be diameter-$k$-critical). We will give a bound of  the number of some edges in $G$,  which is a lower order term compared to $n^2$.

For  $2 \leq i \leq k$, we say that an edge $e$, a pair $\{x,y\}$,  the shortest $(x,y)$-paths   are \emph{$i$-associated} if $d_G(x,y) \leq i$ and $d_{G-e}(x,y) > i$. Similarly, a pair $\{x,y\}$ is \emph{$i$-critical} if there exists some edge $e$ that is  $i$-associated with $\{x,y\}$, and   the shortest $(x,y)$-paths are called  \emph{$i$-critical} paths.
Note that if $i=2$, there exists only one 2-critical $(x,y)$-path; if $i \geq 3$, there may exist more than one $i$-critical $(x,y)$-path of the same length, and these $i$-critical $(x,y)$-paths must contain all the $i$-associated edges.
If not necessary, we simply call a pair $\{x,y\}$ and  the shortest $(x,y)$-path \emph{critical} when there exists an edge $e$ that is $i$-associated with them for some $i$. For any $i$-critical pair $\{x,y\}$,
we choose one of the
$i$-critical $(x,y)$-paths to be $Q_{xy}^i$.
Let $\mathcal{Q}^i(e)$ be the set of   all $i$-critical paths $Q_{xy}^i$ which are $i$-associated with $e$. The \emph{multiplicity} of an edge $e$ is defined as $m(e):=\sum_{i=2}^k|\mathcal{Q}^i(e)|$.

\begin{lem}\label{furedi} For any integer $t$, at most $\frac{k(k+1)}{2t} \binom{n}{2}$ edges have multiplicity at least $t$.
\end{lem}

{\noindent{\bf{Proof.}}}
Since  the $i$-critical path is $i$-associated with at most $i$ edges, we get $\sum_e|\mathcal{Q}^i(e)| \leq i \binom{n}{2}$. Thus $\sum_em(e) =\sum_{i=2}^k\sum_e|\mathcal{Q}^i(e)| \leq \frac{k(k+1)}{2} \binom{n}{2}$, implying that at most $\frac{k(k+1)}{2t} \binom{n}{2}$ edges has multiplicity at least $t$. \bqed

\vskip 2mm

Next, we bound the number of edges whose multiplicity  is less than $t$.
Let $\mathcal{P}_t^i$ ($i \geq 2$) be the set of such $i$-critical paths $Q_{xy}^i$ with length $i $ that  every  edge  incident with $x$ or $y$ is $i$-associated with $Q_{xy}^i$, and has multiplicity less than $t$. For convenience, we call such edges in $Q_{xy}^i$ \emph{t-edges}.
Recall that
a \emph{ $3$-uniform} hypergraph $H$ is a pair $(V(H),E(H))$, where the edge-set $E(H)$ is a collection of 3-element subsets of $V(H)$, and we call the 3-element set in $E(H)$ as \emph{3-edge}. We say a hypergraph $H$ is \emph{linear} if any two distinct 3-edges share at most one vertex.  In a linear 3-uniform hypergraph, three 3-edges form a \emph{triangle} if they form a structure isomorphic to $\{\{1,2,3\},\{3,4,5\},\{5,6,1\}\}$.
Let $RSz(n)$ be the maximum number of 3-edges in a triangle-free, linear 3-uniform hypergraph on $n$ vertices. We now introduce the  famous $(6,3)$ theorem of Ruzsa and Szemer\'edi \cite{RS}, which is important  to bound the number of  $t$-edges in $G$.

\begin{thm}\label{RSZ} (Ruzsa and Szemer\'edi \cite{RS})
$RSz(n)=o(n^2)$.
\end{thm}

 F\"{u}redi \cite{Fu} proved that $|\mathcal{P}_t^2|  \leq 27t*RSz(n)$ for diameter-2-critical graphs.  Loh and Ma \cite{LM} proved that for diameter-3-critical graphs, if  the middle edges of all 3-critical paths in $\mathcal{P}_t^3$ have multiplicity at least $t$, then $|\mathcal{P}_t^3|  \leq 54t*RSz(n)$. We generalize their results in the following lemma for graphs with diameter $k \geq 2$.

\begin{lem}\label{lohma} $|\mathcal{P}_t^i|  \leq 54t*RSz(n)$ for $i \geq 2$. As a result, the number of t-edges  is at most $108kt*RSz(n)$.
\end{lem}

{\noindent\bf{Proof.}}
 Define the 3-uniform hypergraph $H_1$ with $V(H_1)=V(G)$, and $E(H_1)$ is formed by arbitrarily choosing exactly one of  $\{\{x,a_1,y\},$
 $\{x,a_{i-1},y\}\}$ of every  $P:=xa_1\cdots a_{i-1}y $ in $ \mathcal{P}_t^i$, so that $|E(H_1)|=|\mathcal{P}_t^i|$. If $i=2$, then we call one of $\{x,y\}$ the \emph{handle},  $a_1$ the \emph{center}.
 If $i \geq 3$ and we choose $\{x,a_1,y\}$ of $P:=xa_1\cdots a_{i-1}y $, we call $x$ the \emph{handle}, $a_1$ the \emph{center}.

Obviously, any two 3-edges of $H_1$ intersect at most two vertices. We first prove that for any 3-edge $\{x,a,y\}$ with the handle $x$ and  the center $a$, the number of other  3-edges in $H_1$ intersecting $\{x,a,y\}$ in exactly 2 vertices is at most $2t-4$. Observe that any other 3-edge can not contain both $x$ and $y$. Since $m(xa)<t$, there exist at most $t-2$ other   3-edges containing $\{x,a\}$.
If $i=2$, then there also exist at most $t-2$ other   3-edges containing $\{a,y\}$. Now we suppose that $i \geq 3$.
Let $b$ be the neighbor of $y$ in the corresponding $i$-critical path  $P \in \mathcal{P}_{t}^i$ of $\{x,a,y\}$.
Since $by$ is $i$-associated with $P$, $aPy$ is an $(i-1)$-critical path of length $i-1$, and is $(i-1)$-associated with $by$,
implying that any other $i$-critical path in $\mathcal{P}_t^i$ containing $\{a,y\}$ must contain $by$. Since $m(by)<t$, there are   at most $t-2$ other 3-edges in $H_1$ containing $\{a,y\}$. The above arguments imply that the number of 3-edges in $H_1$ intersecting $\{x,a,y\}$ in exactly 2 vertices is at most $2t-4$.


We now construct a linear 3-uniform hypergraph $H_2$ from $H_1$ step by step: for each step, we choose  such 3-edge $\{x,a,y\}$ that there exists  3-edge intersect it in exactly two vertices. Delete all these 3-edges intersecting  $\{x,a,y\}$ in two vertices. At the end of this process, we obtain a linear 3-uniform hypergraph $H_2$ such that $|E(H_2)| \geq |E(H_1)|/{2t}=|\mathcal{P}_t^i|/{2t}$.

A 3-uniform hypergraph $H$ is called \emph{3-partite} if there is a partition $V(H)=V_1\cup V_2\cup V_3$ such that each 3-edge intersects each $V_j$ exactly one vertex for all $j \in \{1,2,3\}$.

\begin{fac}\label{fact} (Erd\H{o}s and Kleitman \cite{EK}) Any $r$-uniform hypergraph $H$ contains a $r$-partite subhypergraph $H'$
such that $|E(H')| \geq \frac{r!}{r^r}|E(H)|$.
\end{fac}

Applying Fact \ref{fact} with $r=3$, we obtain a 3-partite subhypergraph $H_3$ of $H_2$ with parts $V_1,V_2,V_3$, such that $|E(H_3)| \geq \frac{2}{9}|E(H_2)|$.

Without loss of generality, we assume that at least $1/6$ of the 3-edges of $H_3$ have handles in $V_1$ and centers in $V_2$. Then we can choose a subhypergraph $H_4$ of $H_3$ such that $|E(H_4)| \geq \frac{1}{6}|E(H_3)|$ and all handles and centers are  in $V_1$ and $V_2$, respectively.

By the above inequalities of $E(H_j)$ for $j \leq 4$, we derived that $|\mathcal{P}_t^i| \leq 54t*|E(H_4)|$. Note that $H_4$ is a linear 3-uniform 3-partite hypergraph with parts $V_1,V_2,V_3$, where the handles are in $V_1$ and the centers are in $V_2$. By the definition of $RSz(n)$, it remains to prove that $H_4$ has no triangles.
 Suppose there exists a triangle with three 3-edges $T_1,T
_2,T_3$ in $H_4$.
By the definition of a triangle, we get $|V(T_1)\cup V(T
_2)\cup V(T_3)|=6$ and   $|V_j \cap (V(T_1)\cup V(T
_2)\cup V(T_3))|=2$ for $1 \leq j \leq 3$. Let $V_j \cap (V(T_1)\cup V(T
_2)\cup V(T_3))=\{a_j,b_j\}$. Without loss of generality, suppose  $T_1=\{a_1,a_2,a_3\}$, $T_2=\{a_1,b_2,b_3\}$ and $T_3 = \{b_1,b_2,a_3\}$. By the definition of 3-edge, there exist  $i$-critical paths $a_1a_2\cdots a_3$, $a_1b_2\cdots b_3$ and $b_1b_2\cdots a_3$ of length $i$.  But  there is a path $a_1b_2\cdots a_3$ of length $i$ contained in  the subgraph $a_1b_2\cdots b_3 \cup b_1b_2\cdots a_3$,  and it does not pass through $a_1a_2$, a contradiction to the fact that  $a_1a_2$ is $i$-associated with $a_1a_2\cdots a_3 $. \bqed
\vskip 2mm

Let $t := \sqrt{n^2/RSz(n)}$. By Lemmas \ref{furedi} and \ref{lohma},
 the sum of the number of all edges with multiplicity at least $t$, and the number of $t$-edges  is $k^2o(n^2)$.
Let $G_0$ be the subgraph obtained from $G$ by deleting these edges. Then $e(G)=e(G_0)+k^2o(n^2)$.

\section{Proof of Theorem \ref{n^2/2k}}\setcounter{equation}{0}\label{theorem2}

In this section, let $G$ be a  diameter-$k$-critical graph with $k \geq 3$. Since $e(G)=e(G_0)+k^2o(n^2)$,  we are left to prove that $e(G_0)\leq \frac{n^2}{2k}+k^2o(n^2)$. We will give some  lemmas of $G_0$ to prove it.

\begin{lem} \label{eG0-k-c-path}
Any two edges of $E(G_0)$ cannot be $i$-associated with  the same $i$-critical path for any $i\geq 2$.
\end{lem}

\noindent{\bf{Proof}.}
Suppose to the contrary that  there exists an $i$-critical path $P=x\cdots ab\cdots cd \cdots y$, where $ab,cd \in E(G_0)$  are  both $i$-associated with $P$. By the definition of critical paths,   $P_1 = abP cd$  is an $i_1$-critical path of length  $i_1 \leq i$, which is   $i_1$-associated with both $ab$ and $cd$. Since $m(ab)<t$ and $m(cd) <t$, we see that both $ab$ and $cd$ are $t$-edges, a contradiction to $ab,cd \in E(G_0)$. \qed
\vskip 2mm
By using Lemma \ref{eG0-k-c-path}, we have the following lemma.

\begin{lem} \label{k-1}
We can find at least $(k-1)(e(G_0)- \frac{n}{2})$ distinct critical pairs.
\end{lem}

\noindent{\bf{Proof}.}
Note that every edge  is $k$-associated with some $k$-critical paths of length at most $k$.
For any edge $e \in E(G_0)$ that is $k$-associated with $k$-critical path $P_e$ of length $k$, we can find a
subpath of $P_e$ of length $i$ containing $e$, which is also an $i$-critical path of length $i$, and is $i$-associated with $e$ for every  $2 \leq i \leq k$. That is, we find $k-1$ distinct critical paths for each of these edges.

For any edge $e \in E(G_0)$ that is $k$-associated with $k$-critical path $P_e$ of length $k-1$, but not  with a $k$-critical path  of length $k$, obviously $e$ is a 2-critical path. We can find
subpath of $P_e$ of length $i$ containing $e$, which is also an $i$-critical path of length $i$, and is $i$-associated with $e$ for every  $2 \leq i \leq k-1$. That is, we find $k-1$ distinct critical paths for every such edge.

 By Lemma \ref{eG0-k-c-path}, all the above critical paths are distinct. Actually, all  the above critical paths
 have distinct endpoint-pairs. This means that we find $k-1$ distinct  critical pairs for every edge  above.

For any two edges  $ab,cd \in E(G_0)$ that are $k$-associated with $k$-critical paths $P_1,P_2$ of length at most $k-2$, respectively, we denote $P_1=x_1\cdots ab\cdots y_1$. Note that $ab$ is $k$-associated with the $k$-critical pair $\{x_1,y_1\}$. Suppose   $ab$ and $cd$ share a common vertex.
Without loss of generality, we only consider $a=c$.
By Lemma \ref{eG0-k-c-path},  $bad$ is not a 2-critical path, i.e., $bd \in E(G)$ or there exists another vertex $a'$ such that $ba'd$ is a path of length two. By replacing the edge $ab$ with $ adb$ or  $ ada'b$ in $P_1$, we obtain a path  of length at most $k$ connecting $x_1$ and $y_1$, which  contradicts the fact that $ab$ is $k$-associated with $\{x_1,y_1\}$. Hence, all edges  that are $k$-associated with $k$-critical paths  of length at most $k-2$ form a matching, implying that the number of these edges is at most $\frac{n}{2}$. This completes the proof.
 \qed

\vskip 2mm
For any graph $F$, let $Disj(F)=\{\{x,y\}~|~N_{F}(x)\cap N_{F}(y)=\emptyset\}$.  For $G_0$, we have the following lemma.

\begin{lem} \label{disj}
Every critical pair is contained in $Disj(G_0)$.
\end{lem}

\noindent{\bf{Proof.}}
For  any critical pair $\{x,y\}$, if $d_G(x,y) = 1$ or $d_G(x,y) \geq 3$, it is obvious that $N_{G}(x)\cap N_{G}(y)=\emptyset$. Thus, $N_{G_0}(x)\cap N_{G_0}(y)=\emptyset$.  If $d_G(x,y)=2$, we see that  $\{x,y\}$ is  a 2-critical pair, and  denote the 2-critical path as $xay$. By the definition of $G_0$, at most one edge of $xay$  belongs to $G_0$, implying that $N_{G_0}(x)\cap N_{G_0}(y)=\emptyset$. \qed
\vskip 2mm

By Lemmas \ref{k-1} and \ref{disj}, we get $|Disj(G_0)| \geq (k-1)(e(G_0)-\frac{n}{2})$.
Finally, we need the  important lemma from F\"uredi \cite{Fu}.

\begin{lem} \label{e+disj} (F\"uredi \cite{Fu})
For any graph $F$ with $n$ vertices, $e(F)+|Disj(F)| \leq \frac{n^2}{2}$.
\end{lem}

By using Lemma \ref{e+disj} to $G_0$, we get $e(G_0) \leq \frac{n^2}{2k}+ \frac{n}{2}$, completing the proof of Theorem \ref{n^2/2k}. \bqed

\vskip 5mm
\noindent{\bf Remarks.}
There is another interesting problem for diameter-$k$-critical graphs. Cacetta and H\"aggkvist \cite{CH}  conjectured that for any diameter-2-critical graph $G$, $\sum_v d_v^2 \leq n*e(G)$. It is easy to check that if this conjecture is true, then $e(G) \leq \frac{n^2}{4}$,  solving the first half of the Murty-Simon Conjecture.  But Loh and Ma \cite{LM} found an infinite family of diameter-2-critical graphs $G$ with $\sum_v d_v^2 \geq (\frac{10}{9}-o(1))n*e(G)$, where $o(1)$ tends to 0 as $n$ tends to infinity. They also proved that there exist absolute constants $c$ and $N$ such that $\sum_v d_v^2 \leq (\frac{6}{5}-c)n*e(G)$ for any diameter-2-critical graph $G$ with $n \geq N$.
 It is interesting to determine  the smallest factor  of $n*e(G)$ in this inequality. For $k \geq 3$,
 Loh and Ma \cite{LM} also proved that $\sum_v d_v^2 \leq n*e(G)$ for any diameter-$k$-critical graph $G$, and this bound is asymptotically tight.

\vskip 5mm

\noindent{\bf\large Acknowledgements}
\vskip 2mm

This research was supported by
 National Key R\&D Program of China (Grant No. 2023YFA1010202), National Natural Science Foundation of China (Grant No. 12401447),  the Open Project Program of Key Laboratory of Discrete Mathematics with Applications of Ministry of Education, Center for Applied Mathematics of Fujian Province, Key Laboratory of Operations Research and Cybernetics of Fujian Universities (Grant No.J20240901), Fuzhou University.
The author Zhu was supported by NSFC under grant number 12401454,  NSF of Jiangsu under grant number BK20241361, and State-sponsored Postdoctoral Researcher program under number GZB20240976 .

\end{document}